\documentclass[11pt]{amsart2000}%arxiv

\newtheorem{proposition}{Proposition}
\newtheorem{corollary}{Corollary}
\theoremstyle{definition}
\newtheorem{definition}{Definition}

\begin{document}
\setlength{\unitlength}{0.01in}
\linethickness{0.01in}
\begin{center}
\begin{picture}(474,66)(0,0)
\multiput(0,66)(1,0){40}{\line(0,-1){24}}
\multiput(43,65)(1,-1){24}{\line(0,-1){40}}
\multiput(1,39)(1,-1){40}{\line(1,0){24}}
\multiput(70,2)(1,1){24}{\line(0,1){40}}
\multiput(72,0)(1,1){24}{\line(1,0){40}}
\multiput(97,66)(1,0){40}{\line(0,-1){40}}
\put(143,66){\makebox(0,0)[tl]{\footnotesize Proceedings of the Ninth Prague Topological Symposium}}
\put(143,50){\makebox(0,0)[tl]{\footnotesize Contributed papers from the symposium held in}}
\put(143,34){\makebox(0,0)[tl]{\footnotesize Prague, Czech Republic, August 19--25, 2001}}
\end{picture}
\end{center}
\vspace{0.25in}
\setcounter{page}{9}
\title{Some equivalences for Martin's Axiom in asymmetric topology}
\author{Bruce S. Burdick}
\address{Roger Williams University\\
Bristol, Rhode Island 02809}
\email{bburdick@rwu.edu}
\thanks{Bruce S. Burdick,
{\em Some equivalences for Martin's Axiom in asymmetric topology},
Proceedings of the Ninth Prague Topological Symposium (Prague, 2001),
pp.~9--13, Topology Atlas, Toronto, 2002}
\begin{abstract}
We find some statements in the language of asymmetric topology and 
continuous partial orders which are equivalent to the statements 
$\kappa < {\mathfrak m}$ or $\kappa < {\mathfrak p}$. 
\end{abstract}
\maketitle

We think of asymmetric topology as those parts of topology in which the 
specialization order, $x \leq y$ if and only if $x \in c\{y\}$, need not
be symmetric. 
(See \cite{Ko} for some of the motivations.) 

Martin's Axiom has many equivalent statements, consequences, and
variations in the literature which can be stated in topological terms. 
Most of the treatments we have seen so far from set-theoretic
topologists assume that spaces are Hausdorff. In view of recent
interest in asymmetric topology, in which even $T_1$ spaces are a highly
symmetric special case, we give some equivalences for Martin's Axiom
which utilize the terms of this field. 

Our reference for properties related to Martin's Axiom is \cite{Fr}, and
for properties related to continuous lattices we referred to \cite{LM}. 

\begin{definition}
A partially ordered set $(\mathcal{P},\leq)$ is {\it upwards-ccc} if any
uncountable subset of $\mathcal{P}$ must have two distinct members which have
a common upper bound in $\mathcal{P}$. 
The cardinal $\mathfrak m$ is the least cardinal such that there exists a
non-empty upwards-ccc partially ordered set $(\mathcal{P},\leq)$ and a
collection $\{D_{\alpha}|\ \alpha < \mathfrak m \}$ of cofinal subsets of
$\mathcal{P}$ such that no upwards-directed subset of $\mathcal{P}$ meets each
$D_{\alpha}$.
\end{definition} 

It can be shown that $\omega_1 \leq {\mathfrak m} \leq {\mathfrak c}$. 
The Martin's Axiom of the title is the statement $\mathfrak m = \mathfrak c$. 

\begin{definition}
A topological space is {\it ccc} if any uncountable collection of open
sets has two distinct members which are not disjoint. 
A space is {\it locally compact} if every open set contains a compact
neighborhood of each of its points. 
For two sets $U$ and $V$, we say $U$ {\it is compact in} $V$, denoted 
$U \ll V$, if for every open cover of $V$ there is a finite subcollection
which covers $U$. 
A space is {\it core compact} if every open set $O$ is the union of open
sets $O'$ with $O' \ll O$. 
\end{definition}

\begin{definition} 
A subset $C$ of a space is {\it irreducible} if it is closed and the
non-empty open sets in the subspace topology on $C$ form a filterbase. 
A space is {\it sober} if it is $T_0$ and every irreducible subset is the
closure of a point. 
A space is {\it supersober} if it is $T_0$ and every ultrafilter on the
space has either the empty set or the closure of a point as its set of
limits. 
A filter $\mathcal{F}$ of open sets on a space is {\it Scott open} if whenever
$\mathcal{D}$ is a collection of open sets which is directed upward under 
inclusion, $\bigcup \mathcal{D} \in \mathcal{F}$ implies that some element 
of $\mathcal{D}$ is a member of $\mathcal{F}$. 
\end{definition}

We will make use below of the Hofmann-Mislove Theorem, that on a sober 
space the intersection of a Scott open filter of open sets is 
non-empty. 
(See \cite{HM} for the original reference, and see \cite{KP} for a new
proof of the form of the theorem that we have just stated.) 

\begin{definition}
In a partially ordered set $(\mathcal{P},\leq)$, for two elements $x$ and $y$
we say $x$ {\it is way below} $y$, denoted $x \ll y$, if for any directed
subset $D$ of $\mathcal{P}$ with supremum $\sup D$, if $y \leq \sup D$ then
there is some $d \in D$ with $x \leq d$. 
$(\mathcal{P},\leq)$ is {\it upwards-continuous} if for all $y \in \mathcal{P}$ 
the set $\{ x \in \mathcal{P}|\ x \ll y \}$ is a directed set with supremum 
$y$. 
\end{definition}

\begin{definition}
A lattice is {\it distributive} if it satisfies the distributive laws, 
$$a \wedge (b \vee c) = (a \wedge b) \vee (a \wedge c),$$ 
$$a \vee (b \wedge c) = (a \vee b) \wedge (a \vee c).$$
A complete lattice is {\it downwards-ccc} if in any uncountable subset 
there are two distinct members $x$ and $y$ with $x \wedge y \ne 0$, where
0 is the least element. 
An element $x$ in a complete lattice is a {\it non-zero divisor} if $x
\wedge y \ne 0$ for any $y \ne 0$. 
An element $x$ is {\it irreducible} if whenever $x = a \wedge b$ then 
$x = a$ or $x = b$. 
\end{definition}

\begin{proposition}
Let $\kappa$ be a cardinal. 
The following are equivalent. 
\begin{enumerate}
\item 
$\kappa < \mathfrak m$
\item 
In a ccc locally compact sober space the intersection of $\kappa$ or fewer
open dense sets is non-empty. 
\item 
In a ccc core compact space if $\mathcal{D}$ is a collection of $\kappa$ or
fewer open dense sets then there is an irreducible set meeting every
member of $\mathcal{D}$. 
\item 
Let $\mathcal{L}$ be a downwards-ccc upwards-continuous distributive complete
lattice. 
If $D$ is a set of $\kappa$ or fewer non-zero divisors of $\mathcal{L}$ then
there is an irreducible $p \in \mathcal{L}$ such that no element of $D$ is
below $p$. 
\end{enumerate}
\end{proposition}

\begin{proof}
(1) implies (2). 
Given a ccc locally compact sober space $(X,\mathcal{T})$ let 
$\mathcal{P} = \mathcal{T} - \{\phi\}$ and define a partial order on 
$\mathcal{P}$ by $O_1 \leq O_2$ if either $O_1 = O_2$ or $O_2$ is compact 
in $O_1$. 
By local compactness, two open sets $O_1$ and $O_2$ have a common upper
bound in $(\mathcal{P},\leq)$ if and only if $O_1 \cap O_2 \ne \phi$. 
So $(\mathcal{P},\leq)$ is upwards-ccc. 
Given a collection $\{O_{\alpha}|\ \alpha < \kappa \}$ of open dense
subsets of $X$, for each $\alpha$ define 
$D_{\alpha} = \{ O \in \mathcal{P}|\ O \subseteq O_{\alpha} \}$. 
Then each $D_{\alpha}$ is cofinal in $(\mathcal{P},\leq)$, again by local
compactness. 
Since $\kappa < \mathfrak m$ there is a directed set $S \subseteq \mathcal{P}$
that meets each $D_{\alpha}$. 
We may assume that $S$ is closed upwards under $\subseteq$, i.e., 
$S$ is a filter of open sets on $(X,\mathcal{T})$. 

If $\bigcap S \ne \phi$ we are done since an element of $\bigcap S$ must be 
contained in some member of each $D_{\alpha}$ and so is in 
$\bigcap_{\alpha < \kappa} O_{\alpha}$. 
If $S$ has a maximal element under $\leq$ then it is $\bigcap S$ and we 
are done. 
Otherwise $S$ is Scott open, since if $\mathcal{D}$ is a directed 
collection of open sets with $\bigcup \mathcal{D} \in S$ then some element 
$O \in S$ is strictly greater than $\bigcup \mathcal{D}$, which means it 
is compact in $\bigcup \mathcal{D}$, and so it is a subset of some 
element of $\mathcal{D}$. 
By the Hofmann-Mislove Theorem $\bigcap S \ne \phi$. 

(2) implies (4). 
Given a downwards-ccc upwards-continuous distributive complete lattice
$\mathcal{L}$, by Theorem 1.2 in \cite{LM} there exists a locally compact
sober space $(X,\mathcal{T})$ such that $\mathcal{L}$ is order isomorphic to
$\mathcal{T}$ with inclusion. 
We note that this $(X,\mathcal{T})$ must be a ccc space. 
Given $\{ d_{\alpha}|\ \alpha < \kappa \}$, a set of non-zero divisors,
the image of each $d_{\alpha}$ under the isomorphism must be an open dense
subset $O_{\alpha}$ of $X$. 
So there is a point $x \in \bigcap_{\alpha < \kappa} O_{\alpha}$. 
Let $O = X - c\{x\}$. 
Then $O$ corresponds in $\mathcal{L}$ to an irreducible $p$, and since $x$ is
in each $O_{\alpha}$ we can't have $d_{\alpha} \leq p$ for any $\alpha$. 

(4) implies (3). 
Given a ccc core compact space $(X,\mathcal{T})$, by Theorem 1.1 in \cite {LM}
$(\mathcal{T},\subseteq)$ is an upwards-continuous lattice. 
It is clear that it is also downwards-ccc, distributive, and 
complete. 
Given a collection $\{ O_{\alpha}|\ \alpha < \kappa \}$ of open dense
subsets of $X$, we see as above that each $O_{\alpha}$ is a non-zero
divisor. 
So there is an irreducible $O$ in the lattice which does not contain any
$O_{\alpha}$. 
$C = X - O$ is an irreducible set for $(X,\mathcal{T})$, and $C$ meets every
$O_{\alpha}$. 

(3) implies (1). 
Given an upwards-ccc partially ordered set $(\mathcal{P},\leq)$ consider the
partial order topology $\mathcal{T}$ generated by all sets 
$$\uparrow x = \{ y \in \mathcal{P}|\ x \leq y \}$$ 
for $x \in \mathcal{P}$. 
$(\mathcal{P},\mathcal{T})$ is a ccc space. 
It is locally compact, hence core compact. 
Given a collection $\{ C_{\alpha}|\ \alpha < \kappa\}$ of cofinal subsets
of $\mathcal{P}$, each $C_{\alpha}$ is dense, so there is an irreducible set
$C$ meeting each $C_{\alpha}$. 
But in the partial order topology a set is irreducible only if it is
directed, and so we are done.
\end{proof}

We note that each of properties (2) or (3) in Proposition 1 also implies 
the Baire category formulation of $\kappa < \mathfrak m$, to wit, no ccc 
compact Hausdorff space is the union of $\kappa$ or fewer nowhere dense
sets. 
We gave in the proof an argument that (3) implies (1) so that we may now
suggest that (3) is a simultaneous generalization of the partial order and
Baire category forms of $\kappa < \mathfrak m$. 

We also wish to point out that the equivalence of properties (2) and (3)
may be established directly using the sobrification of the space
$(X,\mathcal{T})$. 
The sobrification of $(X,\mathcal{T})$ is the collection of irreducible subsets
of $X$, topologized by the lower Vietoris topology, and $(X,\mathcal{T})$ can
be mapped continuously to its sobrification by sending each $x$ to
$c\{x\}$ (and thus the image of $(X,\mathcal{T})$ is the $T_0$-identification
of $(X,\mathcal{T})$). 
A space is core compact if and only if its sobrification is locally
compact \cite{LM}. 

In view of the Baire category form of $\kappa < \mathfrak m$ mentioned above,
the core compactness in property (3) may be replaced by any stronger
property that is weaker than compact Hausdorff. 
Likewise, locally compact sober in property (2) may be replaced by
anything stronger which is implied by compact Hausdorff, including locally
compact supersober. 
Sober may also be weakened to {\it quasisober}, which is sober without the
$T_0$ assumption, since that is sufficient for the Hofmann-Mislove
Theorem. 

\begin{corollary}
The following are equivalent.
\begin{enumerate} 
\item 
Martin's Axiom. 
\item 
In a ccc locally compact sober space the intersection of fewer than $\mathbb
c$ open dense sets is non-empty. 
\item 
In a ccc core compact space if $\mathcal{D}$ is a collection of fewer than
$\mathfrak c$ open dense sets then there is an irreducible set meeting every
member of $\mathcal{D}$. 
\item 
Let $\mathcal{L}$ be a downwards-ccc upwards-continuous distributive complete
lattice. 
If $D$ is a set of fewer than $\mathfrak c$ non-zero divisors of $\mathcal{L}$
then there is an irreducible $p \in \mathcal{L}$ such that no element of $D$
is below $p$. 
\end{enumerate}
\end{corollary}

There is another cardinal, $\mathfrak p$, for which we can get some similar
results. 

\begin{definition}
In a partially ordered set $(\mathcal{P},\leq)$ a set $R$ is {\it 
upwards-centered} if every finite subset of $R$ has an upper bound in 
$\mathcal{P}$. 
$(\mathcal{P},\leq)$ is {\it upwards}-$\sigma$-{\it centered} if 
$\mathcal{P}$ is the union of countably many centered subsets. 
The cardinal $\mathfrak p$ is the least cardinal such that there exists a 
non-empty upwards-$\sigma$-centered partially ordered set $(\mathcal{P},\leq)$
and a collection $\{D_{\alpha} | \alpha < \mathfrak p\}$ of cofinal subsets of
$\mathcal{P}$ such that no upwards-directed subset of $\mathcal{P}$ meets
each $D_{\alpha}$.\footnote{This is not the original definition. 
if $\mathcal{A}$ is a collection of infinite sets we say that an infinite set
$B$ is a {\it pseudo-intersection} of $\mathcal{A}$ if $B - A$ is finite for
every $A \in \mathcal{A}$. 
Then $\mathfrak p$ is the least cardinal such that there exists a collection
$\mathcal{A}$ of subsets of $\omega$, with the cardinaltiy of
$\mathcal{A}$ equal to $\mathfrak p$, and although every finite subset of
$\mathcal{A}$ has an infinite intersection, $\mathcal{A}$ has no infinite
pseudo-intersection. 
It was Murray Bell \cite{Be} who proved that this definition is equivalent
to the $\sigma$-centered partial order definition given above, and Fremlin 
\cite{Fr} refers on page 25 to this result as ``Bell's Theorem.''}
A space $(X,\mathcal{T})$ is $\sigma$-{\it centered} if 
$\mathcal{T} - \{\phi\}$ with reverse inclusion is
upwards-$\sigma$-centered. 
A complete lattice $\mathcal{L}$ is {\it downwards}-$\sigma$-{\it centered} if
$\mathcal{L} - \{0\}$ with order reversed is upwards-$\sigma$-centered.
\end{definition}

The Baire category formulation of $\kappa < \mathfrak p$ is that no
separable compact Hausdorff space is the union of fewer than $\mathfrak p$
nowhere dense sets. 
It can be shown that $\mathfrak m \leq \mathfrak p \leq \mathfrak c$. 

\begin{proposition} 
Let $\kappa$ be a cardinal. 
The following are equivalent. 
\begin{enumerate}
\item 
$\kappa < \mathfrak p$
\item 
In a separable locally compact sober space the intersection of $\kappa$ or
fewer open dense sets is non-empty. 
\item 
In a separable core compact space if $\mathcal{D}$ is a collection of $\kappa$
or fewer open dense sets then there is an irreducible set meeting every
member of $\mathcal{D}$. 
\item 
In a $\sigma$-centered locally compact sober space the intersection of
$\kappa$ or fewer open dense sets is non-empty. 
\item 
In a $\sigma$-centered core compact space if $\mathcal{D}$ is a collection of
$\kappa$ or fewer open dense sets then there is an irreducible set meeting
every member of $\mathcal{D}$. 
\item 
Let $\mathcal{L}$ be a downwards-$\sigma$-centered upwards-continuous 
distributive complete lattice. 
If $D$ is a set of $\kappa$ or fewer non-zero divisors of $\mathcal{L}$ then
there is an irreducible $p \in \mathcal{L}$ such that no element of $D$ is
below $p$. 
\end{enumerate}
\end{proposition}

\begin{proof}
The equivalence of properties (1), (4), (5), and (6) may be established as
in the proof of Proposition 1. 
Properties (2) and (3) are implied by (4) and (5), respectively, and each
of them implies the Baire category form of property (1).
\end{proof}

We do not know if $\sigma$-centered is equivalent to separable for locally
compact sober spaces, but this is true for locally compact supersober
spaces. 

%\bibliographystyle{amsplain}
%\bibliography{02}
\providecommand{\bysame}{\leavevmode\hbox to3em{\hrulefill}\thinspace}
\providecommand{\MR}{\relax\ifhmode\unskip\space\fi MR }
\providecommand{\MRhref}[2]{%
  \href{http://www.ams.org/mathscinet-getitem?mr=#1}{#2}
}
\providecommand{\href}[2]{#2}

\end{document}